\magnification\magstep1
\baselineskip = 18pt
\def\n{\noindent}
\input mssymb.tex
\null\vskip.1in

\centerline{\bf Ellipsoids of maximal volume in convex bodies}\bigskip
\centerline{Keith Ball}\bigskip
\centerline{Department of Mathematics}
\centerline{Texas A\&M University}
\centerline{College Station, TX \ 77843}\bigskip\medskip

\n {\bf Abstract.} The largest discs contained in a regular tetrahedron lie
in its faces. The proof is closely related to the theorem of Fritz John
characterising ellipsoids of maximal volume contained in convex bodies.
\bigskip

\n {\bf \S 0. Introduction.}

In 1948, Fritz John proved that each convex body in ${\Bbb R}^n$ contains
an unique ellipsoid of maximal volume. Thus, each convex body has an affine
image whose ellipsoid of maximal volume is the Euclidean unit ball,
$B^n_2$. John characterised these affine images with the following theorem.

\proclaim Theorem [J]. The Euclidean ball is the ellipsoid of maximal
volume contained in the convex body $C \subset {\Bbb R}^n$ if and only if
$B^n_2 \subset C$ and, for some $m\ge n$, there are Euclidean unit vectors
$(u_i)^m$, on the boundary of $C$, and positive numbers $(c_i)^m_1$ for
which\medskip
\item{a)} $\sum\limits_i c_iu_i = 0$ and
\item{b)} $\sum\limits_i c_iu_i \otimes u_i = I_n$, the identity on ${\Bbb
R}^n$.$\hfill \square$\medskip

The $u_i$'s of the theorem, are points of contact of the unit sphere
$S^{n-1}$ with the boundary of $C$. The theorem says that weights may be
distributed on the collection of such points so that, a) \ the centre of
mass of the distribution is at the origin and b) \ the inertia tensor of
the distribution is the identity. The first condition shows that the
contact points do not all lie ``on one side'' of the sphere, and the
second, that they do not all lie ``close to proper subspace''.

Condition b) shows that the $u_i$'s behave like an orthonormal basis in
that the inner product can be represented

$$\langle x,y\rangle = \sum_i c_i\langle u_i,x\rangle \langle u_i,y\rangle,
\quad x,y\in {\Bbb R}^n.$$

\n It follows immediately from the equality of the traces of the operators
in b) that

$$\sum_i c_i = n.\eqno (1)$$

At each $u_i$, the supporting hyperplane to $C$ (is unique and) is
perpendicular to $u_i$ (since this is true for the Euclidean ball). Hence,
the set $K = \{x\in {\Bbb R}^n\colon \ \langle x,u_i\rangle \le 1, \ 1 \le
i \le m\}$ contains $C : C$ is sandwiched between $B^n_2$ and $K$. From
this it follows (and this was one of the motivations for John's theorem)
that $C$ is contained in $nB^n_2$, the Euclidean ball of radius $n$. To see
this, suppose $x\in K$ and $\|x\| = r$. Since $-r \le \langle x,u_i\rangle
\le 1$ for every $i$,

$$\eqalign{0 &\le \sum_i c_i(1-\langle x,u_i\rangle) (r + \langle
x,u_i\rangle)\cr
&= r \sum_i c_i + (1-r) \sum_i c_i \langle x,u_i\rangle - \sum_i c_i\langle
x,u_i\rangle^2.}$$

\n Properties a) and b) and (1) show that the latter is

$$rn - \|x\|^2 = rn-r^2.$$

\n Hence $r\le n$.

John's theorem has been used many times in the theory of finite-dimensional
normed spaces. For symmetric convex bodies (the unit balls of normed
spaces) condition a) is redundant and a stronger estimate $\|x\| \le
\sqrt{n}$ can be obtained for $x\in K$. As a consequence, every
$n$-dimensional normed space is isomorphic, with isomorphism constant at
most $\sqrt{n}$, to $n$-dimensional Euclidean space.

It was pointed out to me a few months ago by Prof. A.~Pe\'lcy\'nski that
the
literature does not contain any very explicit proof of the easy (if) part
of John's theorem. The first section of this article contains a quick proof
of this assertion. At about the same time, several people (not quite
independently) asked me for a proof of the intuitively obvious fact that
the largest discs contained in a regular tetrahedron, lie in its faces. The
second section of this article consists of a proof of (the analogue of)
this fact for $k$-dimensional Euclidean balls inside regular
$n$-dimensional simplices. For each $n$ and $k<n$, the regular
$n$-dimensional simplex which circumscribes $B^n_2$, contains a
$k$-dimensional Euclidean ball of radius

$$\sqrt{n(n+1)\over k(k+1)} = r(n,k)\quad {\rm (say),}$$

\n in each of its $k$-dimensional faces. The relationship between Sections~
1 and 2 of this article is elucidated in Section~3 where it is shown that
if $C$ is a convex body in ${\Bbb R}^n$ whose ellipsoid of maximal volume
is $B^n_2$, then $C$ does not contain $k$-dimensional ellipsoids whose
volume is larger than that of a $k$-dimensional ball of radius $r(n,k)$.

The result of Section~3, and, a fortiori, that of Section~2, could
certainly be ``checked by hand''. It is enough to show that the convex hull
of $B^n_2$ and a ``large'' $k$-dimensional ellipsoid, contains an
$n$-dimensional ellipsoid of volume larger than $B^n_2$. But the
calculations involved are messy. The argument presented in Section~3 is a
compelling, if simple, illustration of the usefulness of John's
characterisation.\medskip

\n {\bf \S 1. The `if' part of John's theorem.}

\proclaim Proposition 1. Let $(u_i)^m_1$ be a sequence of unit vectors in
${\Bbb R}^n$ and $(c_i)^m_1$ a sequence of positive numbers satisfying
\medskip
\item{a)} $\sum\limits_i c_iu_i = 0$ and
\item{b)} $\sum\limits_i c_iu_i \otimes u_i = I_n$.\par
\n Then the set $K = \{x\in {\Bbb R}^n\colon \ \langle x,u_i\rangle \le 1,
\ 1 \le i \le m\}$ contains an unique ellipsoid of maximal volume, the
Euclidean unit ball.\medskip

\n {\bf Proof.} Let ${\cal E}$ be the ellipsoid,

$$\left\{x\in {\Bbb R}^n\colon \  \sum^n_1 \alpha^{-2}_j \langle x-y,
v_j\rangle^2 \le 1\right\}$$

\n for some $y \in {\Bbb R}^n$, orthonormal basis $(v_j)^n_1$ and positive
numbers $(\alpha_j)^n_1$. The problem is to show that if ${\cal E} \subset
K$, then $\Pi \alpha_j \le 1$ with equality only if $\alpha_j = 1$ for all
$j$, and $y=0$.

Now, for each $i, \ 1 \le i \le m$, the point

$$x_i = y + \left( \sum^n_{j=1} \alpha^2_j\langle
u_i,v_j\rangle^2\right)^{-{1\over 2}} \sum^n_{j=1} \alpha^2_j \langle
u_i,v_j\rangle v_j$$

\n belongs to ${\cal E}$ and so $\langle u_i,x_i\rangle \le 1$ for each
$i$. Hence

$$\langle u_i,y\rangle + \left( \sum^n_{j=1} \alpha^2_j \langle
u_i,v_j\rangle^2\right)^{1\over 2} \le 1\eqno (2)$$

\n for each $i$. Multiply by $c_i$, sum over $i$ and use the fact that
$\sum\limits_i c_iu_i = 0$ to get

$$\sum_i c_i\left( \sum_j \alpha^2_j \langle u_i,
v_j\rangle^2\right)^{1\over 2} \le \sum_i c_i= n.$$

\n Since $\sum\limits_i c_i\langle u_i,x\rangle^2 = \|x\|^2$ for all $x$,
and the $v_j$'s form an orthonormal basis,

$$\eqalign{\sum_j \alpha_j &= \sum_j\sum_i \alpha_jc_i \langle u_i,v_j
\rangle^2\cr
&= \sum_i c_i \left( \sum_j \alpha_j \langle u_i, v_j\rangle^2\right)\cr
&\le \sum_i c_i\left(\sum_j \alpha^2_j \langle u_i,v_j\rangle^2\right)
^{1\over 2} \left( \sum_j \langle u_i,v_j\rangle^2\right)^{1\over 2}\cr
&= \sum_i c_i\left( \sum_j \alpha^2_j \langle u_i, v_j\rangle^2
\right)^{1\over 2} \le n.}$$

\n By the $AM\backslash GM$ inequality $\Pi \alpha_j \le 1$. There is
equality only if $\alpha_j = 1$ for all $j$ in which case (2) says that

$$\langle u_i,y\rangle + \|u_i\| \le 1 \quad \hbox{for all}\quad i,$$

\n i.e. $\langle u_i,y\rangle \le 0$ for all $i$. Since $\sum c_i\langle
u_i, y\rangle = 0$, this implies that $\langle u_i,y\rangle = 0$ for all
$i$ and so $y = 0$.$\qquad \square$\vfill\eject

\n {\bf \S 2. The simplex.}

\proclaim Proposition 2. Let $T$ be a regular solid simplex in ${\Bbb R}^n$
of internal radius 1. For $1\le k\le n-1$, the largest $k$-dimensional
Euclidean balls in $T$ are those of radius
$$\sqrt{n(n+1)\over k(k+1)}$$
which lie in $k$-dimensional faces of $T$.

The proof of Propositon~2 makes use of the following well-known application
of Caratheodory's theorem.

\proclaim Lemma 3. If $(x_i)^m_1$ is a sequence in ${\Bbb R}^k$ of diameter
at most $d$, there is a point $x\in {\Bbb R}^k$ with
$$
\|x_i-x\| \le r = d\sqrt{k\over 2(k+1)} \quad \hbox{for all}\quad i.
$$
The bound is sharp only if $(x_i)^m_1$ includes some $k+1$ points, all at a
distance $r$ from their average.

\n {\bf Proof.} Let $x$ be the point of ${\Bbb R}^k$ which minimises
$\max\limits_i \|x_i-x\|$ and suppose that this maximum is $s$. Than $x$ is
in the convex hull of those $x_i$'s from which it has distance $s$ since
otherwise there would be a small perturbation, $y$, of $x$ with
$\max\limits_i\| x_i-y\| < s$. By Caratheodory's theorem, $x$ is a convex
combination of some $k+1$ of the $x_i$'s at  distance $s$ from $x\colon$ \
say

$$x = \sum^{k+1}_1 \lambda_ix_i.$$

\n Then

$$\eqalign{s^2 &= \sum^{k+1}_1 \lambda_i \|x_i-x\|^2\cr
&= {1\over 2} \sum_{i,j} \lambda_i\lambda_j \|x_i-x_j\|^2\cr
&= {1\over 2} \sum_{i\ne j} \lambda_i\lambda_j \|x_i-x_j\|^2\cr
&\le {1\over 2} d^2 \sum_{i\ne j} \lambda_i\lambda_j = {1\over 2} d^2
\left( \left(\sum \lambda_i\right)^2 - \sum \lambda_i^2\right)\cr
&\le {d^2k\over 2(k+1)}}$$

\n since $\sum \lambda^2_i \ge {1\over k+1} \big(\sum \lambda_i\big)^2 =
{1\over k+1}$ by the Cauchy-Schwartz inequality.

For the estimate to be sharp, one needs that $\lambda_i  = {1\over k+1}, \
1 \le i\le k+1$, implying the second assertion of the lemma.\medskip$\hfill
\square$

\n {\bf Proof of Proposition 2.} Assume that $T$ is given by

$$T = \{x\in {\Bbb R}^n\colon \ \langle x,u_i\rangle \le 1, \ 1 \le i \le
n+1\}$$

\n for an appropriate sequence $(u_i)^{n+1}_1$ of unit vectors and note
that

$$\eqalignno{&\sum_i u_i = 0&(3)\cr
&\sum_i u_i\otimes u_i = {n+1\over n} I_n.&(4)}$$

\n Let ${\cal E} = \{x \in {\Bbb R}^n\colon \ \sum\limits^k_1 \langle x-y,
v_j\rangle^2 \le r^2, \langle x-y, v_j\rangle = 0, k+1 \le j \le n\}$ be a
$k$-dimensional ball contained in $T$, for some $y\in {\Bbb R}^n, \ r>0$
and orthonormal basis $(v_j)^n_1$. As in the proof of Proposition~1,

$$\langle u_i,y\rangle + r \left(\sum^k_{j=1} \langle u_i, v_j\rangle^2
\right)^{1\over 2} \le 1\quad \hbox{for all}\quad i.$$

\n Let $P$ be the orthogonal projection of ${\Bbb R}^n$ onto ${\rm
span}(v_j)^k_1$. Then for all $i$

$$\langle u_i,y\rangle + r \|Pu_i\| \le 1.\eqno (5)$$

\n Summing over $i$ and using (3),

$$r \sum_i \|Pu_i\| \le n+1$$

\n so it is enough to show that

$$\sum_i \|Pu_i\| \ge \sqrt{k(k+1)(n+1)\over n}\eqno (6)$$

\n for every orthogonal projection of rank $k$.

Now, the set $\{u_i\}^{n+1}_1$ has diameter

$$\sqrt{2(n+1)\over n}$$

\n and hence the set $\{Pu_i\}^{n+1}_1$ has diameter at most this. Since
$\{Pu_i\}^{n+1}_1$ sits in an Euclidean space of dimension $k$, Lemma~3
shows that there is a point $x\in P({\Bbb R}^n)$ with

$$\|Pu_i-x\| \le \sqrt{k(n+1)\over (k+1)n} \quad \hbox{for all}\quad i.
\eqno (7)$$

From identity (4),

$$\sum_i Pu_i\otimes Pu_i = {n+1\over n} P$$

\n and equating traces,

$$\sum_i \|Pu_i\|^2 = {n+1\over n} k.$$

\n Since $\sum\limits_i Pu_i  = 0$,

$$\eqalign{{n+1\over n} k&= \sum_i \langle Pu_i, Pu_i\rangle\cr
&= \sum_i \langle Pu_i, Pu_i-x\rangle\cr
&\le \sum_i \|Pu_i\| \cdot \|Pu_i-x\|\cr
&\le \sqrt{k(n+1)\over (k+1)n} \sum_i \|Pu_i\|,}$$

\n giving the desired inequality (6).

Now, suppose the maximum radius is attained. Then there is equality in (7)
so, by Lemma~3 again, the set $\{Pu_i\}$ includes $k+1$ points all at
distance $\sqrt{k(n+1)\over (k+1)n}$ from their average. But every $k+1$
of the $u_i$'s are at this distance from their average. So $P$ is an
isometry on the affine hull of some $k+1, \ u_i$'s: \ i.e. the range of $P$
is parallel to this affine hull. This implies that ${\cal E}$ lies in a
$k$-dimensional subspace parallel to some $k$-dimensional face of $T$. This
fact determines all the numbers $\|Pu_i\|$ and hence the numbers $\langle
u_i,y\rangle$ since there is equality in (5) for all $i$. These numbers
determine $y$.$\hfill \square$\bigskip

\n {\bf \S 3. The general case.}

\proclaim Proposition 4. Let $(u_i)^m_1, \ (c_i)^m_1$ and $K$ be as in
Proposition~1. If ${\cal E}$ is a $k$-dimensional ellipsoid in $K$ then the
($k$-dimensional) volume of ${\cal E}$ is no larger than that of a
$k$-dimensional ball of radius $\sqrt{n(n+1)\over k(k+1)}$.

This proposition cannot be proved using the argument of Proposition~2 as it
stands, since for a general sequence $(u_i)^m_1$ satisfying the hypotheses
and orthogonal projection $P$ of rank $k, \ \sum\limits_i c_i\|Pu_i\|$ may
be
as small as $k$. This complicates the argument somewhat: \ Proposition~2 is
isolated because it has a simpler proof. For the proof of Proposition~4,
Lemma~3 is replaced by an easier observation.

\proclaim Lemma 5. If $(x_i)^m_1$ is a sequence of vectors with $\sum x_i =
0$ and $(u_i)^m_1,$ a sequence of unit vectors in some Euclidean space then

$$\left(\sum_i \langle x_i,u_i\rangle\right)^2 \le \sum_{i,j} \|x_i\| \cdot
\|x_j\| (1-\langle u_i,u_j\rangle).$$

\n {\bf Proof.} By homogeneity, it may be assumed that $\sum\limits_i
\|x_i\| = 1$. Set $\lambda_i = \|x_i\|, \ 1 \le i\le m$ and $u =
\sum\limits_i \lambda_iu_i$. Then,

$$\eqalignno{\left(\sum_i \langle x_i,u_i\rangle\right)^2 &= \left(\sum
\langle x_i,u_i-u\rangle \right)^2\cr
&\le \left(\sum \|x_i\| \cdot \|u_i-u\|\right)^2\cr
&= \left(\sum \lambda_i \|u_i-u\|\right)^2\cr
&\le \sum \lambda_i \|u_i-u\|^2 = 1 - \|u\|^2\cr
&= \sum_{i,j} \lambda_i\lambda_j (1-\langle u_i,u_j\rangle).&\square}$$
\medskip

\n {\bf Proof of Proposition 4.} Let ${\cal E}$ be the ellipsoid

$$\eqalign{\bigg\{ x\in {\Bbb R}^n\colon \ &\sum^k_1 \alpha^{-2}_j \langle
x-y, v_j\rangle^2 \le 1,\cr
&\langle x-y, v_j\rangle = 0, \ k+1\le j\le n\bigg\}}$$

\n for some $y\in {\Bbb R}^n$, orthonormal basis $(v_j)^n_1$ and positive
numbers $(\alpha_j)^k_1$. The problem is to show that

$$\left(\prod^k_1 \alpha_j\right)^{1\over k} \le \sqrt{n(n+1)\over
k(k+1)}.$$

\n It certainly suffices to show that

$$\sum^k_1 \alpha_j \le \sqrt{kn(n+1)\over k+1}.$$

As in the proof of Proposition 1,

$$\langle u_i,y\rangle + \left( \sum^k_{j=1} \alpha^2_j \langle u_i,v_j
\rangle^2 \right)^{1\over 2} \le 1$$

\n for every $i$. Define $T\colon \ {\Bbb R}^n\to {\Bbb R}^n$ by

$$Tx = \sum^k_{j=1} \alpha_j \langle x,v_j\rangle v_j.$$

\n Then

$$\langle u_i,y\rangle + \|Tu_i\| \le 1\quad \hbox{for every}\quad i.\eqno
(8)$$

\n Also

$$\sum^k_1 \alpha^2_j = \sum_{i,j} \alpha^2_j c_i \langle u_i,v_j\rangle^2
= \sum_i c_i\|Tu_i\|^2\eqno (9)$$

\n and

$$\sum^k_1 \alpha_j = \sum_{i,j} \alpha_j c_i\langle u_i,v_j
\rangle^2
= \sum_i c_i \langle Tu_i,u_i\rangle.\eqno (10)$$

\n The proof divides into two parts, the first of which effectively handles
the $u_i$'s which are far apart (as if the body were symmetric) while the
second, and more complicated part, handles the $u_i$'s which are close
together (as if the body were  a simplex).

Since $\|Tu_i\| \le 1 - \langle u_i,y\rangle$ for each $i$ (by (8))

$$\eqalign{\sum_i c_i\| Tu_i\|^2 &\le \sum_i c_i (1-\langle
u_i,y\rangle)^2\cr
&= \sum_i c_i - \left\langle \sum_i c_iu_i,y\right\rangle + \sum_i c_i
\langle u_i,y\rangle^2\cr
&= n+\|y\|^2.}$$

\n So by (9)

$${1\over k} \left(\sum_i \alpha_j\right)^2 \le \sum_j \alpha^2_j \le n +
\|y\|^2.\eqno (11)$$

\n On the other hand, set $x_i = c_iTu_i, \ 1 \le i\le m$ and observe that
$\sum x_i = T\big(\sum c_iu_i\big) = 0$. Then Lemma~5 and (10) show that

$$\eqalign{\left(\sum^k_1 \alpha_j\right)^2 &= \left(\sum_i \langle
x_i,u_i\rangle\right)^2\cr
&\le \sum_{i,j} \|x_i\|\cdot \|x_j\| (1-\langle u_i,u_j\rangle)\cr
&= \sum_{i,j} c_ic_j (1-\langle u_i,u_j\rangle) \|Tu_i\|\cdot \|Tu_j\|.}$$

\n Since $c_ic_j(1-\langle u_i,u_j\rangle) \ge 0$ for all $i$ and $j$, (8)
can be applied again to give

$$\left(\sum \alpha_j\right)^2 \le \sum_{i,j} c_ic_j (1-\langle
u_i,u_j\rangle ) (1-\langle u_i,y)\rangle (1-\langle u_j,y\rangle).$$

\n Expanding this product and using the fact that $\sum\limits_i c_iu_i =
0$ one obtains

$$\left(\sum \alpha_j\right)^2 \le \left(\sum c_i\right)^2 - \sum_{i,j}
c_ic_j \langle u_i,u_j\rangle \langle u_i,y\rangle \langle u_j,y\rangle.$$

\n Two applications of the identity
$\sum\limits_i c_i\langle u_i,x\rangle \langle u_i,y\rangle = \langle
x,y\rangle$ show that

$$\left(\sum \alpha_j\right)^2 \le n^2 - \|y\|^2.\eqno (12)$$

\n Finally, this inequality may be added to (11) to give

$$\left(1 + {1\over k}\right) \left(\sum \alpha_j\right)^2 \le n^2+n,$$

\n and hence

$$\left( \sum \alpha_j\right)^2 \le {kn(n+1)\over k+1}$$

\n as required.$\hfill \square$\medskip

\n {\bf Remarks.} For $k=1$, Proposition~4 states that if $C$ is a convex
body whose ellipsoid of maximal volume is $B^n_2$ then ${\rm diam}(C)\le
\sqrt{2n(n+1)}$. This fact can be proved more simply: \ if $x,y\in C$ then
$\|x\|, \|y\| \le n$, as explained in the introduction, and (with the usual
notation)

$$\eqalign{0 &\le \sum_i c_i (1-\langle u_i,x\rangle )(1-\langle
u_i,y\rangle)\cr
&= n + \langle x,y\rangle}$$

\n so that

$$\eqalign{\|x-y\|^2 &= \|x\|^2 + \|y\|^2 - 2\langle x,y\rangle\cr
&\le 2n^2 + 2n = 2n(n+1).}$$

The fact that $\|x\| \le n$ for $x\in C$ could be deduced from inequality
(12) of the above proof, with ``$k=0$''.\bigskip\medskip

\centerline{\bf References}

\item{[J]} F. John, Extremum problems with inequalities as subsidiary
conditions, Courant Anniversary Volume, Interscience, New York, 1948,
187-204.

\end